\documentclass[12pt]{amsart}

\usepackage{amsthm}
\usepackage{amssymb}
\usepackage{amsmath}
\numberwithin{equation}{section}
\newcommand{\bea}{\begin{eqnarray}}
\newcommand{\eea}{\end{eqnarray}}
\newcommand{\be}{\begin{eqnarray*}}
\newcommand{\ee}{\end{eqnarray*}}
\newtheorem{theorem}{Theorem}[section]
\newtheorem{lemma}{Lemma}[section]
\newtheorem{corollary}{Corollary}[section]
\newtheorem{definition}{Definition}[section]
\newtheorem{proposition}{Proposition}[section]
\newtheorem{example}{Example}[section]

\begin{document}

\title[Coincidence Symmetry Groups]{Structures of Coincidence Symmetry Groups}
\author[Yi Ming Zou]{Yi Ming Zou}
\address{Department of Mathematical Sciences, University of Wisconsin, Milwaukee, WI 53201, USA} \email{ymzou@uwm.edu}
\maketitle

\begin{abstract}
The structure of the coincidence symmetry group of an arbitrary $n$-dimensional lattice in the $n$-dimensional Euclidean space is considered by describing a set of generators. Particular attention is given to the coincidence isometry subgroup (the subgroup formed by those coincidence symmetries which are elements of the orthogonal group). Conditions under which the coincidence isometry group can be generated by reflections defined by vectors of the lattice will be discussed, and an algorithm to decompose an arbitrary element of the coincidence isometry group in terms of reflections defined by vectors of the lattice will be given.
\end{abstract}

\section{Introduction}
\par
The mathematical theory of coincidence site lattice (CSL) can be used to describe certain phenomena that arise in the physics of interfaces and grain boundaries (for a more detailed background in CSL theory, we refer the readers to the references, especially Baake [1997], Bollmann [1970], and Grimmer [1973]). Because of the success of the models for crystalline interfaces based on the properties of CSL and related lattices (Brandon {\it et al} [1964]; Bollmann [1970]; Warrington \& Bufalini [1971]; Grimmer [1973], [1976]), the focus of the CSL theory has been mostly on the coincidence of two lattices of the same dimensions (the coincidence of two lattices of different dimensions can be easily reduced to the same dimension case). M. A. Fortes [1983] developed a matrix theory of CSL by using the normal form of an integer matrix. In the first paper of [1983], Fortes gave a crystallographic interpretation of the invariant set of an integer matrix, and applied it to solve the degree of coincidence problem of two lattices in arbitrary dimensions. In the subsequence, the theory was extended to include displacement shift complete (DSC) lattices, and a method to calculate bases for these lattices via some special factorizations of the related matrices was provided. Duneau {\it et al.} [1992] further developed the matrix theory of CSL also by using the normal form of an integer matrix, and gave a method to decompose the corresponding matrix into associated shear transformations. Pleasants {\it et al.} [1996] used number theory to solve the planar coincidences for $N$-fold symmetry. Baake [1997] used the factorization properties of certain number fields to solve the coincidence problem for dimensions up to 4. Recently, Arag\'{o}n {\it et al.} [2001] and Rodriguez {\it et al.} [2005] developed a different approach to coincidence isometry theory by using geometric algebra (Clifford algebra) as a tool. From the work of these literatures, problems on the structures of the coincidence symmetry group of a given lattice can be formulated. In this paper, we consider the structure of the coincidence isometry group of a lattice in $\mathbb{R}^{n}$.
\par
Let $L$ be a lattice with basis $(a_{1}, \ldots, a_{n})$, let $V$ be the $n$-dimensional real vector space with the same basis, let $\mathcal{A}$ be a linear transformation of $V$, and let $A$ be the matrix of $\mathcal{A}$ under the basis $(a_{1}, \ldots, a_{n})$. We call $\mathcal{A}$ a {\it coincidence symmetry} if $\mathcal{A}$ is an automorphism of $V$ and $L\cap \mathcal{A}L$ is a sublattice of $L$ with finite index. If $\mathcal{A}$ is a coincidence symmetry of $L$, we call $A$ a coincidence matrix of $L$, or abusing language, we also call $A$ a coincidence symmetry. It is known (see section 2 below) that $A$ is a coincidence symmetry if and only if $A$ is a rational matrix. The set of all coincidence symmetry (or the set of all $n\times n$ coincidence matrices) of $L$ forms a group under the multiplication defined by composition (or the multiplication of matrices). If $L$ is a lattice of the Euclidean space $\mathbb{R}^{n}$, then one can consider the isometries of $\mathbb{R}^{n}$ which are coincidence symmetries of $L$. In this case, one has the {\it coincidence isometry} subgroup formed by all the coincidence isometries (Baake [1997]). We analyze the structures of these groups by considering the decomposition of a matrix from both geometric and algebraic view points. Baake [1997] (see also Pleasants {\it et al.} [1996]) uses the factorization of numbers to reduce a symmetry to irreducible ones, while the approach developed by Arag\'{o}n {\it et al.} [2001] and Rodriguez {\it et al.} [2005] relies on the decomposition of a matrix into product of coincidence reflections. The results in Arag\'{o}n {\it et al.} [2001] stated that if the matrix is a product of coincidence reflections, then the corresponding symmetry is a coincidence isometry. In Rodriguez {\it et al.} [2005], it was conjectured that any coincidence isometry of the lattice spanned by the canonical basis of $\mathbb{R}^{n}$ is a product of coincidence reflections. We shall prove a theorem which includes this conjecture as a special case, and use the theorem to describe the coincidence isometry group. 
\par
In section 2, we briefly recall the relevant definitions and some known results. In section 3, we prove a theorem about coincidence isometry groups of lattices $L$ in $\mathbb{R}^{n}$, and apply it to describe the structure of the coincidence isometry group. Examples will be given in section 4.  

\section{Notation and definitions}

The set of real numbers is denoted by $\mathbb{R}$, the set of real $n\times n$ matrices is denoted by $M_{n}(\mathbb{R})$, and the set of all non-singular $n\times n$ real matrices is denoted by $GL_{n}(\mathbb{R})$. Notation for matrices over the rational numbers $\mathbb{Q}$ and the integers $\mathbb{Z}$ are defined similarly. For example, $GL_{n}(\mathbb{Z})$ denotes the set of all invertible $n\times n$ integer matrices, so
\be
GL_{n}(\mathbb{Z})=\{\text{$n\times n$ integer matrices $A$ with $\det A=\pm 1$}\}.
\ee
We also consider the above sets of non-singular matrices as linear transformations. For example, we also regard $GL_{n}(\mathbb{R})$ as the set of all non-singular linear transformations of $\mathbb{R}^{n}$. If we do regard them as linear transformations, we will specify the basis which relates the transformations to their matrices.
\par
By an $n$-dimensional lattice $L$ with basis $(a_{1}, \ldots, a_{n})$, we mean the free abelian group $\oplus_{i=1}^{n}\mathbb{Z}a_{i}$. With the basis $(a_{1}, \ldots, a_{n})$, we can always define a standard inner product on the $n$-dimensional real vector space $\oplus_{i=1}^{n}\mathbb{R}a_{i}$ by requiring $(a_{1}, \ldots, a_{n})$ to be an orthonormal basis. This defines an isometry between the usual $n$-dimensional Euclidean space $\mathbb{R}^{n}$ and $\oplus_{i=1}^{n}\mathbb{R}a_{i}$. However, usually we need to consider a lattice in the $n$-dimensional Euclidean space $\mathbb{R}^{n}$ with canonical basis $(e_{1}, \ldots, e_{n})$. In this case, we assume the lattice to be also $n$-dimensional, since if the lattice has dimension $m < n$, then we can always consider the $m$-dimensional subspace of $\mathbb{R}^{n}$ that contains the lattice of interest. Thus, a lattice $L\subset \mathbb{R}^{n}$ is given by an $n\times n$ non-singular matrix $A$ and a basis of the lattice is  
\bea
(a_{1}, \ldots, a_{n})= (e_{1}, \ldots, e_{n})A.
\eea
We call the matrix $A$ the structure matrix of $L$, and use the notation $L_{A}$ if we want specify the fact that the lattice $L$ is given by the matrix $A$.
\par
We adopt the definition that a sublattice $L^{\prime}\subset L$ is a subgroup $L^{\prime}$ of finite index in the abelian group $L$. In the usual notation, this is $[L:L^{\prime}]<\infty$. The CSL theory concerns the problems which arise when the intersection $L_{1}\cap L_{2}$ of two lattices happens to be a sublattice of both lattices $L_{1}$ and $L_{2}$. If this is the case, we say that $L_{1}$ and $L_{2}$ are {\it commensurate lattices}. 
\par
Suppose that $L_{i}$ is given by the structure matrix $A_{i}, i=1,2$, let the basis of $L_{i}$ be $\mathbf{B}_{i}$. Then 
\be
\mathbf{B}_{i}=(e_{1}, \ldots, e_{n})A_{i}, \quad i=1,2. 
\ee
\begin{theorem} [Grimmer] The lattices $L_{1}$ and $L_{2}$ are commensurate if and only if $A_{2}^{-1}A_{1}$ is a rational matrix.
\end{theorem}
\begin{proof}
Let $L^{\prime} = L_{1}\cap L_{2}$ and let $\mathbf{B}^{\prime}$ be a basis of $L^{\prime}$. Then there are integer matrices $N_{i}$ ($i=1,2$) such that 
\be
\mathbf{B_{1}}N_{1} = \mathbf{B}^{\prime} = \mathbf{B_{2}}N_{2}.
\ee
Under the assumption that $L_{1}$ and $L_{2}$ are commensurate, i.e. $[L_{i}:L^{\prime}]<\infty$ ($i=1,2$), the matrices $N_{i}$ are non-singular, thus from $A_{1}N_{1} = A_{2}N_{2}$, we obtain $A_{2}^{-1}A_{1} = N_{2}N_{1}^{-1}$, implies that $A_{2}^{-1}A_{1}$ is a rational matrix. Conversely, if $A_{2}^{-1}A_{1}$ is a rational matrix, then there exists an integer $m>0$ such that $mA_{2}^{-1}A_{1}$ is an integer matrix, say $A$. Then from $mA_{1}=A_{2}A$, we have $m\mathbf{B}_{1}=\mathbf{B}_{2}A$. Hence $mL_{1}\subset L^{\prime}$, which implies that $[L_{1}:L^{\prime}]\le m^{n}$. Symmetrically, we also have $[L_{2}:L^{\prime}]<\infty$. Therefore $L_{1}$ and $L_{2}$ are commensurate. 
\end{proof} 
\par
Grimmer's theorem immediately implies the following:
\begin{corollary}
Let $L$ be a lattice with basis $(a_{1},\ldots, a_{n})$, and let $A$ be an $n\times n$ non-singular real matrix. Then the lattice with basis $(a_{1},\ldots, a_{n})A$ and the lattice $L$ are commensurate if and only if $A$ is a rational matrix.
\end{corollary}
\par
However, if we view the matrix $A$ in the above corollary as the matrix of a linear transformation, then we need to specify under which basis this matrix is given. In Corollary 2.1, the matrix is given by using the basis $(a_{1},\ldots, a_{n})$. Let us consider a lattice $L$ in $\mathbb{R}^{n}$ with the structure matrix $A$. Let $\mathcal T$ be a linear transformation of $\mathbb{R}^{n}$, and let $T$ be the matrix of $\mathcal T$ under the canonical basis $(e_{1},\ldots, e_{n})$. Then the structure matrix of the lattice $\mathcal{T}(L)$ (the image of $L$ under the transformation $\mathcal T$) is $TA$. Then by Theorem 2.1, the lattice $\mathcal{T}(L)$ and the lattice $L$ are commensurate if and only if $A^{-1}TA$ is rational. This leads to the following definition:  
\begin{definition} Let $L_{A}\subset\mathbb{R}^{n}$ be a lattice with the structure matrix $A$. We call the group $AGL_{n}(\mathbb{Q})A^{-1}$ the coincidence symmetry group (CSG) of $L_{A}$.
\end{definition}  
\par
The isometries of $\mathbb{R}^{n}$ (with the standard inner product $(,)$) which provide commensurate lattices to a lattice $L\subset\mathbb{R}^{n}$ are of special interest (cf. Baake [1997], Arag\'{o}n {\it et al.} [2001], and Rodriguez {\it et al.} [2005]). Let $O(n)$ be the set of orthogonal transformations of $\mathbb{R}^{n}$. The concept of coincidence isometry group was defined in Baake [1997] with the notation $OC(L)$, i.e.,
\be
OC(L)=\{ \mathcal{R}\in O(n): [L: L\cap \mathcal{R}L]<\infty\}.
\ee
For our purpose, we need a definition in terms of matrices under the canonical basis of $\mathbb{R}^{n}$. Let
\be
O_{n}(\mathbb{R})=\{A\in M_{n}(\mathbb{R}): A^{t}A = I\}.
\ee
That is, $O_{n}(\mathbb{R})$ is the set of $n\times n$ orthogonal real matrices.
Suppose that $\mathcal R\in O(n)$ and $[L:L\cap \mathcal{R}(L)]<\infty$. Let $R$ be the matrix of $\mathcal R$ under the canonical basis, then $R\in O_{n}(\mathbb{R})$. From the discussion preceding Definition 2.1, we conclude that the matrix $A^{-1}RA$ is rational. Thus we give the following definition:
\par 
\begin{definition} Let $L_{A}\subset\mathbb{R}^{n}$ be a lattice with the structure matrix $A$. We call the group $O_{n}(\mathbb{R}^{n})\cap (AGL_{n}(\mathbb{Q})A^{-1})$ the coincidence isometry group (CIG) of $L_{A}$.
\end{definition}
\par
Thus, the CIG of $L$ is just the group $OC(L)$, and we will use both terms for our convenience. 
\begin{example}
If $L= \mathbb{Z}^{n}$, then $A = I$, and $OC(L)= O_{n}(\mathbb{Q}):= O_{n}(\mathbb{R}^{n})\cap GL_{n}(\mathbb{Q})$. We call the elements of $O_{n}(\mathbb{Q})$ rational orthogonal matrices.
\end{example}
In the next section we will analyze the structure of the coincidence isometry group of an arbitrary lattice $L$ in $\mathbb{R}^{n}$.
\par
\section
{Decomposition of elements of CIG into reflections}
\par
The decomposition of an element of the CIG of a lattice $L\subset \mathbb{R}^{n}$ is central in the Clifford algebra approach to the coincidence site lattice problem developed in Arag\'{o}n {\it et al.} [2001] and Rodriguez {\it et al.} [2005]. It was conjectured (and proved for the planar lattices) in Rodriguez {\it et al.} [2005] that any coincidence isometry of the canonical lattice $\mathbb{Z}^{n}$ of $\mathbb{R}^{n}$ can be decomposed as a product of {\it coincidence reflections} (reflections that belong to the coincidence isometry group of $L$). Note that for the lattice $\mathbb{Z}^{n}$, the corresponding CIG is $O_{n}(\mathbb{Q})$.  Here, we prove a more general theorem which includes the lattice $\mathbb{Z}^{n}$ as a special case. It should be pointed out that although the Cartan-Dieudonn\'{e} theorem (Porteous [1995, Ch. 5]) says that any orthogonal $n\times n$ real matrix can be decomposed into a product of at most $n$ reflections, it is clear that a statement of coincidence isometries of certain lattices can be decomposed into product of coincidence reflections is not a direct consequence of the Cartan-Dieudonn\'{e} theorem (cf. Example 4.2 below).
\par
\begin{theorem}
Let $L\subset \mathbb{R}^{n}$ be a lattice such that the reflection defined by an arbitrary nonzero vector of $L$ is a coincidence isometry of $L$. Then any coincidence isometry of $L$ can be decomposed as a product of at most $n$ reflections defined by the vectors in $L$.
\end{theorem}
\par
\begin{proof}
Let the structure matrix of $L$ be $A$. Then 
\be
\mathbf{B} = (b_{1}, \ldots, b_{n}) := (e_{1}, \ldots, e_{n})A
\ee
is a basis of $L$. Let $\mathcal R\in O(n)$ be a coincidence isometry of $L$. We use induction on $n$ to prove the theorem. It is clear that the theorem holds for $n=1$. Assume that it holds for all $k$ such that $1\le k<n$ and consider the case $n$. We consider two cases: $\mathcal{R}(b_{1})=b_{1}$ or $\mathcal{R}(b_{1})\ne b_{1}$, separately.
\par
In the first case, let
\be
V = \{x\in \mathbb{R}^{n}: (x,b_{1})=0\}.
\ee
Then $V$ is an $n-1$-dimensional subspace of $\mathbb{R}^{n}$, and $V$ is invariant under $\mathcal{R}$, i.e. $\mathcal{R}(V) = V$. Thus, $\mathcal R$ restricts to an orthogonal transformation $\mathcal{R}^{\prime}$ of the $n-1$-dimensional Euclidean subspace $V$. Compare the orthogonal projection $\mathcal{P} : \mathbb{R}^{n}\longrightarrow V$ defined by $b_{1}$:
\bea
\mathcal{P}(x)= x-\frac{(x,b_{1})}{(b_{1},b_{1})}b_{1}, \quad \text{$\forall x\in \mathbb{R}^{n}$},
\eea
with the reflection of $\mathbb{R}^{n}$ defined by $b_{1}$:
\bea
\mathcal{R}_{b_{1}}(x)= x-\frac{2(x,b_{1})}{(b_{1},b_{1})}b_{1}, \quad \text{$\forall x\in \mathbb{R}^{n}$},
\eea
we can see that under the assumption of the theorem, for each $b_{i}$ ($1<i\le n$), there exists an integer $m_{i}>0$ such that $m_{i}\mathcal{P}(b_{i})\in L$. Let $m=m_{2}\cdots m_{n}$, then $m\mathcal{P}(L)\subset L$. Hence, $\mathcal{R}^{\prime}$ is a coincidence isometry of the $n-1$-dimensional lattice $\mathcal{P}(L)$ (with basis $(\mathcal{P}(b_{2}), \ldots, \mathcal{P}(b_{n}))$) which satisfies the condition of the theorem. Therefore, by induction assumption, $\mathcal{R}^{\prime}$ is a product of $j$ reflections defined by some vectors $y_{1}, \ldots, y_{j}\in\mathcal{P}(L)$ such that $1\le j\le n-1$. Let $x_{i}=my_{i}, 1\le i\le j$. Then all $x_{i}\in L$. Let the reflection of $\mathbb{R}^{n}$ defined by $x_{i}$ be $\mathcal{R}_{i}$, then $\mathcal{R} =\mathcal{R}_{1}\cdots\mathcal{R}_{j}$. Hence the theorem is proved in this case.
\par
In the second case, $\mathcal{R}(b_{1})\ne b_{1}$, thus $a:=\mathcal{R}(b_{1})-b_{1}\ne 0$. Let $\mathcal{R}_{a}$ be the reflection defined by the vector $a$. Since $\mathcal{R}$ is a coincidence isometry of $L$, there exists an integer $t>0$ such that $ta\in L$. However, $\mathcal{R}_{a}=\mathcal{R}_{ta}$, so $\mathcal{R}_{a}$ can be viewed as a reflection defined by a vector in $L$. Consider the coincidence isometry $\mathcal{R}_{a}\mathcal{R}$ of $L$. Note that $(b_{1},b_{1})=(\mathcal{R}(b_{1}), \mathcal{R}(b_{1}))$, we have (it can also be seen easily via a geometric diagram)
\be
\mathcal{R}_{a}\mathcal{R}(b_{1})&=& \mathcal{R}(b_{1})-\frac{2(\mathcal{R}(b_{1}),a)}{(a,a)}a\\
  &=& \mathcal{R}(b_{1})-\frac{2(\mathcal{R}(b_{1}),\mathcal{R}(b_{1})-b_{1})}{(\mathcal{R}(b_{1})-b_{1},\mathcal{R}(b_{1})-b_{1})}(\mathcal{R}(b_{1})-b_{1})\\
  &=& b_{1}.
\ee
Thus by the first case, $\mathcal{R}_{a}\mathcal{R}$ is a product of at most $n-1$ reflections defined by some vectors of $L$, say $\mathcal{R}_{a}\mathcal{R} = \mathcal{R}_{1}\cdots\mathcal{R}_{j}$ with $1\le j\le n-1$. Then since $\mathcal{R}_{a}^{2}=I$, we conclude that $\mathcal{R} = \mathcal{R}_{a}\mathcal{R}_{1}\cdots\mathcal{R}_{j}$ is a product of at most $n$ reflections defined by vectors of $L$. This completes the proof of the theorem.
\end{proof}
\par
Note that the proof of Theorem 3.1 gives a practical way to actually decompose a coincidence isometry into a product of coincidence reflections. We will give an example in section 4.
\par
It turns out that the condition in Theorem 3.1 is sufficient for any application purpose for which the computations involve only rational numbers. The following theorem gives a necessary and sufficient condition for a lattice to satisfy the condition in Theorem 3.1.
\par
\begin{theorem} Let $L\subset \mathbb{R}^{n}$ be a lattice with structure matrix $A =(a_{ij})$, and let $a_{i}, 1\le i\le n$, be the column vectors of $A$. Then every nonzero vector of $L$ defines a coincidence reflection of $L$ if and only if the ratios:
\bea
\frac{(a_{j},a_{i})}{(a_{k},a_{k})}, \quad 1\le i,j,k\le n,
\eea
are all rational.
\end{theorem}
\begin{proof}
If every nonzero vector of $L$ defines a coincidence reflection of $L$, then in particular, every $a_{i}$ ($1\le i\le n$) defines a coincidence reflection of $L$. Let $\mathcal{R}_{i}$ be the reflection defined by $a_{i}$. Then since
\bea
\mathcal{R}_{i}(a_{j}) = a_{j}-\frac{2(a_{j},a_{i})}{(a_{i},a_{i})}a_{i}, \quad \forall j,
\eea
we must have 
\bea
\frac{(a_{j},a_{i})}{(a_{i},a_{i})}\in \mathbb{Q}, \quad 1\le i,j\le n.
\eea
If $(a_{i},a_{k})\ne 0$, then
\be
\frac{(a_{i},a_{i})}{(a_{k},a_{k})}=\frac{(a_{i},a_{k})}{(a_{k},a_{k})}\frac{(a_{i},a_{i})}{(a_{i},a_{k})}
\ee 
is a product of two rational numbers, hence is rational. If $(a_{i},a_{k}) = 0$, consider the reflection $\mathcal{R}_{c}$ defined by $c = a_{i}-a_{k}$. By assumption, $\mathcal{R}_{c}$ is a coincidence reflection of $L$. Thus from 
\be
\mathcal{R}_{c}(a_{i}) = a_{i}-\frac{2(a_{i},c)}{(c,c)}c,
\ee
we have 
\be
\frac{(a_{i},c)}{(c,c)} = \frac{(a_{i},a_{i})}{(a_{i},a_{i})+(a_{k},a_{k})}
                =\frac{1}{1+\frac{(a_{k},a_{k})}{(a_{i},a_{i})}}
\ee
is rational. Hence we also have 
\be
\frac{(a_{i},a_{i})}{(a_{k},a_{k})}\in \mathbb{Q}.
\ee 
Together with (3.5), this proves (3.3).
\par
Conversely, if (3.3) holds, let $x = AX\in L$ be a nonzero vector, where $X = (x_{1}, \ldots, x_{n})^{t}\in \mathbb{Z}^{n}$ is a column vector. Then for any $1\le i\le n$, 
\be
\frac{(a_{i},x)}{(x,x)} = \frac{\sum_{j=1}^{n}x_{j}(a_{i},a_{j})}{\sum_{s,t=1}^{n}x_{s}x_{t}(a_{s},a_{t})}
=\frac{\sum_{j=1}^{n}x_{j}\frac{(a_{i},a_{j})}{(a_{i},a_{i})}}{\sum_{s,t=1}^{n}x_{s}x_{t}\frac{(a_{s},a_{t})}{(a_{i},a_{i})}}
\ee 
is rational. It follows that the reflection defined by $x$ is a coincidence isometry of $L$. This completes the proof of the theorem.
\end{proof}
\par
A useful consequence of Theorem 3.2 is the following:
\begin{corollary} Let $L\subset \mathbb{R}^{n}$ be a lattice with the structure matrix $A$. If $A^{t}A$ is a rational matrix, then every nonzero vector of $L$ defines a coincidence reflection of $L$, and hence every coincidence isometry of $L$ can be decomposed into a product of at most $n$ coincidence reflections defined by the vectors of $L$.
\end{corollary}
\begin{proof}
Keep the notation of Theorem 3.2. Under the assumption that $A^{t}A$ is rational, all $(a_{i},a_{j}), 1\le i,j\le n,$ are rational, hence condition (3.3) holds. 
\end{proof}
\par
A special case of Corollary 3.1 is when the matrix $A$ is rational.
\begin{corollary} If $A$ is rational, then every nonzero vector of $L$ defines a coincidence reflection of $L$, and hence every coincidence isometry of $L$ can be decomposed into a product of at most $n$ coincidence reflections defined by the vectors of $L$.
\end{corollary} 
\par
The decomposition of a coincidence isometry of the lattice $L=\mathbb{Z}^{n}$ into a product of coincidence reflections is just a special case of Corollary 3.2. 
\par
By Theorem 3.1 and Theorem 3.2, we immediately obtain the following:
\begin{theorem} If the structure matrix $A$ of a lattice $L\subset \mathbb{R}^{n}$ satisfies condition (3.3), then $OC(L)$ is generated by the reflections defined by the nonzero vectors of $L$. 
\end{theorem}
\par
As an application we have:
\begin{theorem} For $n>1$, $OC(\mathbb{Z}^{n})$ is infinitely generated.
\end{theorem}
To prove Theorem 3.4, we need the following fact about the rational numbers:
\begin{lemma}
Let $S$ be a finite subset of the rational numbers $\mathbb{Q}$, and let $P$ be the set of all the prime integers that show up in the denominators of the reduced forms of the elements of $S$. If only addition, subtraction, and multiplication are allowed, then $S$ can not produce rational numbers whose denominators of the reduced forms contain prime factors not in $P$.
\end{lemma}
Now we are ready to prove Theorem 3.4.
\begin{proof}
Assume that $n>1$, and let $G = OC(\mathbb{Z}^{n})$. Under the assumption of the theorem, every nonzero vector $v\in \mathbb{Z}^{n}$ generates an element $\mathcal{R}_{v}\in G$. By Corollary 2.1, the matrix of $\mathcal{R}_{v}$ under the canonical basis $(e_{1},\ldots, e_{n})$ is a rational matrix. Since the inverse of an orthogonal matrix is its transpose, $G$ is generated as a group by the rational matrices defined by the reflections of the nonzero vectors of $\mathbb{Z}^{n}$ involving only addition, subtraction, and multiplication of rational numbers. If $G$ is finitely generated, then there is a finite subset $S$ of $G$ whose elements are rational matrices that generates $G$. Let $P$ be the set of all the prime integers which show up in the denominators of the reduced forms of the rational numbers involved in the elements of $S$. To prove the theorem, by Lemma 3.1, we only need to show that there is a nonzero vector $v\in \mathbb{Z}^{n}$ such that the matrix of $\mathcal{R}_{v}$ under the canonical basis involves rational numbers whose reduced forms contain prime factors in the denominators which are not in $P$.
\par
We consider vectors of the form
\be
v=e_{1}+ye_{2}, \quad y\in \mathbb{Z},
\ee
and consider the fraction that shows up in
\bea
\mathcal{R}_{v}(e_{1})=e_{1}-\frac{2}{1+y^{2}}(e_{1}+ye_{2}).
\eea
Suppose $p$ is the largest element in $P$. If we let $y=p_{1}\cdots p_{r}$ be the product of the first $r$ primes $\le p$, then all the prime factors of the denominator of the fraction in (3.6) are not in $P$. This completes the proof of Theorem 3.4.
\end{proof}
\par
It should be pointed out that a detailed analysis the group of $OC(\mathbb{Z}^{2})$ is contained in Baake [1997].
\par
\section{Examples}
\par
We consider two examples in this section. In the first example, we show how to use the procedure in the proof of Theorem 3.1 to decompose a coincidence isometry into a product of coincidence reflections. In the second example, we consider a special type of lattices in $\mathbb{R}^{2}$ and determine their coincidence isometry groups.
\par
{\bf Example 4.1.} Let $L\subset \mathbb{R}^{2}$ be the rhombic lattice defined by the matrix
\be
A = \left(\begin{array}{cc} 1 & \frac{2}{3} \\ 0 & \frac{\sqrt{5}}{3} \end{array}\right).
\ee
Let
\be
R = \left(\begin{array}{cc} -\frac{19}{21} & -\frac{4\sqrt{5}}{21}\\ \frac{4\sqrt{5}}{21} & -\frac{19}{21} \end{array}\right).
\ee
Then $R$ is an orthogonal matrix and 
\be
A^{-1}RA = 
\left(\begin{array}{cc} -\frac{9}{7} & -\frac{4}{7}\\ \frac{4}{7} & -\frac{11}{21} \end{array}\right).
\ee
Thus $R$ is a coincidence isometry of the lattice $L$ (the coincidence index is 21). Denote the column vectors of $A$ by $a_{1}, a_{2}$. Let
\be
b_{1} = R(a_{1})-a_{1}= \left(\begin{array}{cc} \frac{-40}{21}\\ \frac{4\sqrt{5}}{21}\end{array}\right).
\ee
Then the matrix of the reflection $\mathcal{R}_{b_{1}}$ under the canonical basis is
\be
R_{1}=\left(\begin{array}{cc} -\frac{19}{21} & \frac{4\sqrt{5}}{21}\\ \frac{4\sqrt{5}}{21} & \frac{19}{21} \end{array}\right),
\ee
and 
\bea
R_{1}R = \left(\begin{array}{cc} 1 & 0  \\ 0 & -1\end{array}\right).
\eea
Let
\be
b_{2} = 2a_{1}-3a_{2} = \left(\begin{array}{cc} 0\\ -\sqrt{5}\end{array}\right).
\ee
Then $b_{2}$ is a scalar multiple of the projection of $a_{2}$ with respect to the orthogonal projection defined by $a_{1}$. The matrix $R_{2}$ of the reflection defined by $b_{2}$ under the canonical basis is the matrix on the right hand side of (4.1) and $R =R_{1}R_{2}$.
\par
\medskip
{\bf Example 4.2.} Let $L\subset \mathbb{R}^{2}$ be a lattice with the structure matrix
\bea
A = \left(\begin{array}{cc} a & 1\\ 0 & b \end{array}\right),
\eea
where $a$ and $b$ are arbitrary positive real numbers. Let the column vectors of $A$ be $a_{1}, a_{2}$. For this matrix, condition (3.3) is equivalent to 
\be 
a, \frac{a}{1+b^{2}}\in \mathbb{Q}\Longleftrightarrow a, b^{2}\in \mathbb{Q}.
\ee
If this is the case, $OC(L)$ is generated by the reflections defined by the nonzero vector of $L$.
\par
If $a\notin \mathbb{Q}$, but $\frac{a}{1+b^{2}}\in \mathbb{Q}$, then $b^{2}\notin\mathbb{Q}$. To find the condition for a reflection to be a coincidence reflection, we only need to consider vectors of the form $v=xe_{1}+e_{2}, x\in \mathbb{R}$ (these vectors need not be in $L$). Consider
\bea
\frac{(v,a_{1})}{(v,v)}v &=& \frac{ax}{1+x^{2}}v=\frac{x(bx-1)}{b(1+x^{2})}a_{1}+\frac{ax}{b(1+x^{2})}a_{2},\\
\frac{(v,a_{2})}{(v,v)}v &=& \frac{x+b}{1+x^{2}}v=\frac{(x+b)(bx-1)}{ab(1+x^{2})}a_{1}+\frac{x+b}{b(1+x^{2})}a_{2}.\nonumber
\eea
If at least one of $x+b$ and $bx-1$ is 0, then the fractions involved in (4.3) are all rational numbers and the reflection defined by $v$ is a coincidence reflection of $L$. In the first case, the vector $v$ is orthogonal to $a_{2}$; and in the second case, the vector $v$ is parallel to $a_{2}$. Assume that both $x+b$ and $bx-1$ are nonzero, and suppose that $v$ defines a coincidence reflection of $L$. Then the second equation in (4.3) implies that $x\ne 0$. Furthermore, (4.3) implies that
\bea
\frac{ax}{x+b},\quad\frac{bx-1}{a}\in\mathbb{Q}.
\eea
Since $\frac{a}{1+b^{2}}\in\mathbb{Q}$, (4.4) implies that
\bea
\frac{x(1+b^{2})}{x+b},\quad\frac{bx-1}{1+b^{2}}\in\mathbb{Q}.
\eea
The second condition in (4.5) implies that there exists a $q\in\mathbb{Q}$ such that
\bea
x=\frac{q(1+b^{2})+1}{b}.
\eea
Substitute (4.6) into the first condition of (4.5), we have $b^{2}\in\mathbb{Q}$, which contradicts our assumption. Thus, the only coincidence reflections are defined by a vector which is parallel to $a_{2}$ or a vector which is perpendicular to $a_{2}$.
\par
Similarly, we can discuss the case that $a\in\mathbb{Q}$ but $b^{2}\notin\mathbb{Q}$ and the case that both $a,\frac{a}{1+b^{2}}\notin\mathbb{Q}$. In the first case, the only coincidence reflections are defined by $a_{1}$ or a nonzero vector which is orthogonal to $a_{1}$. In the second case, there is no coincidence reflection for $L$. To determine the group $OC(L)$, it remains to consider rotations. If 
\be
R = \left(\begin{array}{cc} cos\theta & -sin\theta\\ sin\theta & cos\theta \end{array}\right)
\ee 
is a coincidence isometry of $L$, then
\be
R(a_{1}) &=& \left(\begin{array}{cc} a cos\theta\\ a sin\theta \end{array}\right)= xa_{1}+ya_{2}=\left(\begin{array}{cc} ax + y\\ by \end{array}\right),\\
R(a_{2}) &=& \left(\begin{array}{cc} cos\theta-b sin\theta\\ sin\theta+b cos\theta \end{array}\right)= x^{\prime}a_{1}+y^{\prime}a_{2}=\left(\begin{array}{cc} ax^{\prime} + y^{\prime}\\ b y^{\prime} \end{array}\right),
\ee
for some $x,x^{\prime},y,y^{\prime}\in \mathbb{Q}$. In particular, we have
\bea
b y=a sin\theta, \quad a x^{\prime}
=-(b+\frac{1}{b})sin\theta.
\eea
If $sin\theta\ne 0$, (4.7) implies that 
\bea
\frac{b^{2}+1}{a^{2}}=-\frac{x^{\prime}}{y}\in \mathbb{Q}.
\eea
However, if one of $a$ and $\frac{a}{1+b^{2}}$ is rational and the other one is irrational, then (4.8) does not hold. If both are irrational, then 
\be
(ax+y)^{2}+(by)^{2}=a^{2}
\ee
together with (4.8) will also lead to a contradiction. Therefore $sin\theta =0$ and $R=\pm I$.  
\par
\medskip
Summarize, we have
\begin{proposition}
Suppose the structure matrix of $L\subset\mathbb{R}^{2}$ is given by (4.2). 
\par
(1) If $a, b^{2}\in\mathbb{Q}$, then $OC(L)$ is generated by the reflections defined by the nonzero vectors of $L$.
\par
(2) If $a\in\mathbb{Q}$ but $b^{2}\notin\mathbb{Q}$, then $OC(L)=\{\pm I, \pm R_{a_{1}}\}\cong \mathbb{Z}_{2}^{2}$. 
\par
(3) If $a\notin\mathbb{Q}$ but $\frac{a}{1+b^{2}}\in\mathbb{Q}$, then $OC(L)=\{\pm I, \pm R_{a_{2}}\}\cong \mathbb{Z}_{2}^{2}$. 
\par
(4) If $a, \frac{a}{1+b^{2}}\notin\mathbb{Q}$, then $OC(L)=\{\pm I\}\cong\mathbb{Z}_{2}$.
\end{proposition} 
\par
It should be pointed out that there are many ways to decompose an orthogonal matrix into products of reflections. To see this, we just need to note that the identity matrix is the product of any reflection with itself. It can be seen (say, by considering planar lattices and rotations) that for any integer $m>0$, there are orthogonal transformations $\mathcal R$ of $\mathbb{R}^{n}$ such that $\mathcal{R}^{k}$ are not coincidence isometry of the canonical lattice $\mathbb{Z}^{n}$ for all $1\le k<m$, but $\mathcal{R}^{m}$ is a coincidence isometry of $\mathbb{Z}^{n}$. The same is true for reflections, i.e., there are reflections $\mathcal{R}_{i}$ ($1\le i\le m$) such that any partial product of the $\mathcal{R}_{i}$'s is not a coincidence isometry, but the product $\mathcal{R}_{1}\cdots \mathcal{R}_{m}$ is a coincidence isometry.     
\par
\medskip

\end{document}